\documentclass[11pt]{amsart} 
\usepackage{amssymb,amsmath,latexsym,enumerate,graphicx,bbm,mathptmx,microtype,cite,xcolor}
\usepackage{caption}

\newtheorem{theorem}{Theorem} 
\newtheorem*{lrc}{Lonely Runner Conjecture} 

\newtheorem{corollary}[theorem]{Corollary}
\newtheorem{proposition}[theorem]{Proposition}
\newtheorem{conjecture}[theorem]{Conjecture}

\newtheorem{question}[theorem]{Question}

\renewcommand\emptyset{\varnothing}
\renewcommand\phi{\varphi}

\newcommand\conv{\operatorname{conv}}

\DeclareMathOperator*{\lin}{lin}

\newcommand\ZZ{\mathbb{Z}}
\newcommand\RR{\mathbb{R}}

\newcommand\PP{\mathbb{P}}

\newcommand\cC{\mathcal{C}}

\newcommand\cK{\mathcal{K}}

\newcommand\cP{\mathcal{P}}
\newcommand\cQ{\mathcal{Q}}

\newcommand\cZ{\mathcal{Z}}

\newcommand\be{\mathbf{e}}
\newcommand\bm{\mathbf{m}}
\newcommand\bn{\mathbf{n}}

\newcommand\bt{\mathbf{t}}

\newcommand\bv{\mathbf{v}}
\newcommand\bw{\mathbf{w}}
\newcommand\bx{\mathbf{x}}

\newcommand\bzero{\mathbf{0}}
\newcommand\bone{\mathbf{1}}

\begin{document}

\title{Lonely Runner Polyhedra}

\author{Matthias Beck}
\author{Serkan Ho\c{s}ten}
\address{Department of Mathematics\\
         San Francisco State University\\
         San Francisco, CA 94132\\
         U.S.A.}
\email{[mattbeck,serkan]@sfsu.edu}

\author{Matthias Schymura}
\address{Institute for Mathematics\\
         \'{E}cole Polytechnique F\'{e}d\'{e}rale de Lausanne\\
         CH-1015 Lausanne\\
         Switzerland}
\email{matthias.schymura@epfl.ch}


\begin{abstract}
We study the \emph{Lonely Runner Conjecture}, conceived by J\"org M.~Wills in the 1960's: Given positive
integers $n_1, n_2, \dots, n_k$, there exists a positive real number $t$ such that for all $1 \le j \le
k$ the distance of $t \, n_j$ to the nearest integer is at least $\frac{ 1 }{ k+1 }$. Continuing a
view-obstruction approach by Cusick and recent work by Henze and Malikiosis, our goal is to promote a
polyhedral \emph{ansatz} to the Lonely Runner Conjecture. Our results include geometric proofs of some
folklore results that are only implicit in the existing literature, a new family of affirmative
instances defined by the parities of the speeds, and geometrically motivated conjectures whose resolution
would shed further light on the Lonely Runner Conjecture.
\end{abstract}

\keywords{Lonely Runner Conjecture, view-obstruction problems, lattice points in polyhedral cones.}

\subjclass[2010]{Primary 52C07; Secondary 11J71.}

\date{4 March 2019}
 
\thanks{We thank Christian Haase, Martin Henk, Romanos Malikiosis, Daria Schymura, J\"org M.~Wills, Kevin Woods, and an
anonymous referee for helpful conversations.
MS was partially supported by the Swiss National Science Foundation (SNSF) within the project \emph{Convexity, geometry of numbers, and the complexity of integer programming (Nr.~163071)}}

\maketitle


\section{Introduction}

We study the following conjecture raised by J\"org M.~Wills in the 1960's~\cite{willslonelyrunner}.

\begin{lrc}
Given pairwise distinct integers $n_0, n_1, \dots, n_k$, for each $0 \leq i \leq k$ there exists a real number~$t$ such that for all $0 \le j \le k$, $i \ne j$, the distance of $t \, ( n_i - n_j )$ to the nearest integer is at least $\frac{ 1 }{ k+1 }$.
\end{lrc}

Wills originally formulated this conjecture for \emph{real} numbers $n_0, n_1, \dots, n_k$, but it can be relaxed to the rational and thus integral case~\cite{bohmanholzmankleitman,henzemalikiosis}.
The lower bound $\frac{ 1 }{ k+1 }$ is best possible, as the case $n_j = j$ for $0 \le j \le k$ and a classic result of Dirichlet on Diophantine approximation (see, e.g., \cite{casseldiophantine}) show.
The name \emph{Lonely Runner Conjecture}, introduced by Goddyn in~\cite{bieniagoddynetal}, stems from the charming model of $k+1$ runners going at different constant speeds around a
circular track of length~1 (having started at the same place and time); the conjecture says that each of them will at some point have distance at least $\frac{ 1 }{ k+1 }$ to the other runners.
For more on the Lonely Runner Conjecture's history, proofs for $k \le 6$, and connections to Diophantine approximation, view-obstruction problems, and graph theory, see
\cite{betkewills, 
cusicklonelyrunner3,
cusickpomerance, 
bohmanholzmankleitman, 
barajasserra, 
renault, 
taolonelyrunner}. 

A simple observation leads to a more convenient formulation of the problem:
The distance of any two runners at any given time depends only on their relative speeds.
So we may pick a fixed runner, say the one with speed $n_0$, reduce the speed of every runner by $n_0$ and consider only the loneliness of the first runner that is now stagnant.

\begin{lrc}
Given pairwise distinct positive integers $n_1, n_2, \dots, n_k$, there exists a real number $t$ such that for all $1 \le j \le k$ the distance of $t \, n_j$ to the nearest integer is at least $\frac{ 1 }{ k+1 }$.
\end{lrc}

A speed vector $\bn \in \ZZ_{ >0 }^k$ that satisfies the Lonely Runner Conjecture is called a \emph{lonely runner instance}.
Our goal is to derive novel families of lonely runner instances, using a polyhedral-geometric model.
In the next section, we introduce this model by defining the lonely runner polyhedron $\cP(\bn)$.
It turns out to be closely related  to the zonotopes that where constructed in~\cite{henzemalikiosis}.
We illustrate the utility of a polyhedral \emph{ansatz} in Section~\ref{sect:geometricProofs} by providing geometric proofs of some
folklore results that are only implicit in the existing literature, and by obtaining a new family of lonely runner instances in Theorem~\ref{thm:oddlrc} defined by the parities of the speeds.
In Sections~\ref{sect:ProjArgs} and~\ref{sect:crosssectArgs}, we use suitable projections and cross sections of the lonely
runner polyhedron to obtain families of lonely runner instances that are independent of the fastest runner
(Theorems~\ref{thm:mainProjection} and~\ref{thm:mainCrossSection}).
We close in Section~\ref{sec:musings} by discussing geometrically motivated conjectures whose resolution
would shed further light on the Lonely Runner Conjecture.

\section{A Polyhedral Model for Lonely Runners}

Our starting point  is a view-obstruction problem due to Cusick~\cite{cusicklonelyrunner3} which,
based on the second formulation above, is easily seen to be equivalent to the Lonely Runner Conjecture.
It states that for every $\bn \in \ZZ_{ >0 }^k$, the line $\RR \, \bn$ in direction $\bn$ and passing through the origin, intersects the $k$-dimensional cube
\[
 \cC(\bm) \ := \ \bm + \left[\tfrac{1}{k+1},\tfrac{k}{k+1}\right]^k \ = \ \left\{ \bx \in \RR^k : \, m_j + \tfrac{ 1 }{ k+1 } \ \le \ x_j \ \le \ m_j + \tfrac{ k }{ k+1 } \ \text{ for } 1 \le j \le k \right\}
\]
for some $\bm \in \ZZ_{ \ge 0 }^k$. 
Equivalently, the point $\bn$ belongs to the nonnegative span~$\cK (\bm)$ of $\cC(\bm)$.

The set $\cK (\bm)$ is a \emph{polyhedral cone}, that is, a set of the form $\{ \sum_{ j=1 }^n \lambda_j \, \bw_j : \, \lambda_j \ge 0 \}$ for some $\bw_1, \dots, \bw_n \in \RR^k$.
In our case, $\cK (\bm)$ 
is spanned by all vectors of the form
\begin{align}
  (k+1) \, \bm \ + \text{ a vector consisting of $k$'s and $1$'s, } \label{eqn:coneGenerators}
\end{align}
but not all of these are extreme rays.
From basic notions of polyhedral geometry (see, e.g.,~\cite{ziegler}) one obtains 
\begin{align}
  \cK (\bm) \ &= \ \left\{ \bx \in \RR^k : \, \bigl( (k+1) \, m_i + 1 \bigr) \, x_j \ \le \ \bigl( (k+1) \, m_j + k \bigr) \, x_i \ \text{ for } 1 \le i, j \le k \right\} \nonumber \\
              &= \ \left\{ \bx \in \RR^k : \, \frac{ (k+1) \, m_j + 1 }{ (k+1) \, m_i + k } \ \le \ \frac{ x_j }{ x_i } \ \le \ \frac{ (k+1) \, m_j + k }{ (k+1) \, m_i + 1 } \ \text{ for } 1 \le i < j \le k \right\} \label{eq:lrcfractions1} \\
              &= \ \left\{ \bx \in \RR^k : \, \frac{ 1 }{ (k+1) x_j } - \frac{ k }{ (k+1) x_i } \ \le \ \frac{ m_i }{ x_i } - \frac{ m_j }{
x_j } \ \le \ \frac{ k }{ (k+1) x_j } - \frac{ 1 }{ (k+1) x_i } \ \text{ for } 1 \le i < j \le k \right\} . \label{eq:lrcfractions}
\end{align}
The last formulation motivates the definition of the polyhedron
\begin{align}
  \cP(\bn)
  \ := \ &\left\{ \bx \in \RR^k : \, \frac{ 1 }{ (k+1) n_j } - \frac{ k }{ (k+1) n_i } \ \le \ \frac{ x_i }{ n_i } - \frac{ x_j }{ n_j } \ \le
\ \frac{ k }{ (k+1) n_j } - \frac{ 1 }{ (k+1) n_i } \ \text{ for } 1 \le i < j \le k \right\} \nonumber \\
     = \ &\left\{ \bx \in \RR^k : \, \frac{ n_i - k \, n_j }{ k+1 } \ \le \ n_j \, x_i - n_i \, x_j \ \le \ \frac{ k \, n_i - n_j }{ k+1 } \
\text{ for } 1 \le i < j \le k \right\} . \label{eq:pndef}
\end{align}
By construction, the polyhedron $\cP(\bn)$ consists of all points $\bm \in \RR^k$ such that $\bn \in \cK(\bm)$.
Based on the description~\eqref{eqn:coneGenerators} of the generators of $\cK(\bm)$, we get the equivalences
\begin{align*}
\bn \in \cK(\bm) \ \Longleftrightarrow \ &\ \exists \ r_\bv \geq 0 \ \text{ such that } \ \bn = \sum_{\bv \in \{1,k\}^k} r_\bv \, ((k+1) \, m + \bv) \\
\ \Longleftrightarrow \ &\ \exists \ r_\bv \geq 0 \ \text{ such that } \ \bm = \frac{1}{ (k+1) \, \sum_{\bv} r_\bv } \, \bn \ - \ \frac{1}{ (k+1) \, \sum_{\bv} r_\bv } \, \sum_{\bv} r_\bv \, \bv \\
\ \Longleftrightarrow \ &\ \bm \in \RR \, \bn - \frac{1}{k+1} \, \conv\left\{ \bv : \bv \in \{1,k\}^k \right\} \\
\ \Longleftrightarrow \ &\ \bm \in \RR \, \bn - \left[ \tfrac{1}{k+1},\tfrac{k}{k+1} \right]^k .
\end{align*}
This gives the convenient and useful description
\begin{align}
\cP(\bn) \ = \ \RR \, \bn - \left[\tfrac{1}{k+1},\tfrac{k}{k+1}\right]^k . \label{eqnDescriptionPn}
\end{align}
It also shows that the lonely runner polyhedron $\cP(\bn)$ is closely connected to the zonotopes constructed in~\cite[Section~2.3]{henzemalikiosis}.
In fact, up to a linear transformation that maps the projected lattice $\ZZ^k \mid \bn^\perp$ to $\ZZ^{k-1}$, the zonotopes
in~\cite{henzemalikiosis} are of the form
\[
\cZ(\bn) \ = \ \left[\tfrac{1}{k+1},\tfrac{k}{k+1}\right]^k \ \Big| \ \bn^\perp .
\]
Therefore,
\[
\cP(\bn) \mid \bn^\perp \ = \ \cP(\bn) \cap \bn^\perp \ = \ -\cZ(\bn),\ \text{ or equivalently } \ \cP(\bn) \ = \ -\cZ(\bn) + \RR \, \bn .
\]
Summarizing the previous observations, we can reformulate the Lonely Runner Conjecture geometrically.
The equivalence (a) $\Longleftrightarrow$ (c) was derived already by Chen~\cite[Lemma~1]{chen1991onaconjII}, yet not in a polyhedral context.

\begin{proposition}\label{prop:geometricLRC}
Let $\bn \in \ZZ_{ >0 }^k$.
The following are equivalent:
\begin{enumerate}[{\rm (a)}]
 \item $\bn$ is a lonely runner instance;
 \item there exists an $\bm \in \ZZ^k_{\geq 0}$ such that $\bn \in \cK(\bm)$;
 \item $\cP(\bn) \cap \ZZ^k \ne \emptyset$;
 \item $\cZ(\bn) \cap \left( \ZZ^k \mid \bn^\perp \right) \ne \emptyset$.
\end{enumerate}
\end{proposition}

Thus there are two basic ways to prove that a given $\bn$ is a lonely runner instance.
Namely, one can directly construct an $\bm \in
\ZZ_{ \ge 0 }^k$ such that $\bn \in \cK(\bm)$---equivalently, $\bm \in \cP(\bn)$---, or one can indirectly prove that $\cP(\bn) \cap \ZZ^k$ is nonempty.
We will encounter examples for either of these approaches in the sequel.

\section{Selected Geometric Proofs}\label{sect:geometricProofs}

The geometric viewpoint outlined in the last section yields many presumably folklore results on classes of lonely runner instances.
In this section, we exemplify this in some selected settings and we provide new information with Proposition~\ref{prop:monotonicity} and Theorem~\ref{thm:oddlrc} as well.

We start with the illustrative case of two non-stationary runners, that is, $\bn \in \ZZ^2_{>0}$.
Here, the lonely runner polyhedron reduces to the infinite strip
\[
\cP(\bn) \ = \ \left\{ \bx \in \RR^2 : \, n_1 - 2 \, n_2 \leq 3 \, n_2 \, x_1 - 3 \, n_1 \, x_2 \leq 2 \, n_1 - n_2\right\} .
\]
Since we may assume that $\gcd(n_1,n_2) = 1$, we can invoke B\'{e}zout's Lemma and express every multiple of three as $3 \, n_2 \, x_1 - 3 \, n_1 \, x_2$, for some integers $x_1$ and $x_2$.
The set $\{n_1 - 2 \, n_2, \ldots, 2 \, n_1 - n_2\}$ contains $2 \, n_1 - n_2 - (n_1 - 2 \, n_2) + 1 = n_1 + n_2 + 1 \geq 3$ elements, and thus always a multiple of three.
Thus, $\cP(\bn) \cap \ZZ^2 \ne \emptyset$.

The following is well known among experts and often used to foster case distinctions (see, e.g.,~\cite[Equation~(1.2)]{dubickas2011the}).

\begin{proposition}\label{prop:n1lenk}
Suppose $n_1$ is the largest and $n_k$ the smallest coordinate of $\bn \in \ZZ^k_{> 0}$.
Then, $n_1 \le k \, n_k$ if and only if $\bzero \in \cP(\bn)$.
In particular, if $n_1 \le k \, n_k$ then $\bn$ is a lonely runner instance.
\end{proposition}
\begin{proof}
In view of the inequalities in~\eqref{eq:pndef}, $\bzero \in \cP(\bn)$ if and only if $n_i - k \, n_j \leq 0 \leq k \, n_i - n_j$, for all $1 \leq i < j \leq k$.
This means that for each pair $(i,j)$ we have $n_i \leq k \, n_j$ and $n_j \leq k \, n_i$, one of which is redundant.
Our choice of the labeling of the coordinates of $\bn$ implies that this is in fact equivalent to $n_1 \leq k \, n_k$.
\end{proof}

Note that the supposedly extreme vector $\bn = (1, 2, \dots, k)$ mentioned in the introduction satisfies the condition of Proposition~\ref{prop:n1lenk}, and indeed, this vector lies on the boundary of $\cK(\bzero)$.

It is reasonable to expect that a speed vector whose coordinates form a nonincreasing sequence should be contained in a cone $\cK(\bm)$ corresponding to a lattice point $\bm \in \ZZ^k_{\geq0}$ whose coordinates are also nonincreasing.
It turns out that this is in fact necessary.

\begin{proposition}\label{prop:monotonicity}
Suppose $n_1 \geq \dots \geq n_k$.
If $\bm \in \cP(\bn) \cap \ZZ^k_{\geq 0}$, then $m_1 \geq \dots \geq m_k$. 
\end{proposition}
\begin{proof}
Suppose $m_{j+1} \geq m_j + 1$ for some $1 \leq j \leq k - 1$.
Then the left hand side of the defining inequality in~\eqref{eq:pndef} for the pair $(j,j+1)$ implies that
\[
( m_j (k+1) + k ) \, n_{j+1} \geq (m_{j+1} (k+1) + 1) \, n_j \geq (m_j (k+1) + k) \, n_j + 2 \, n_j ,
\]
which contradicts the assumption $n_j \geq n_{j+1}$, since $m_j \geq 0$.
\end{proof}

A likewise simple argument reveals that hard instances of the Lonely Runner Problem are those that contain multiples of every sufficiently small integer.
In the setting of chromatic numbers of distance graphs, such a statement can be found in the work of Eggleton, Erd\H{o}s \& Skilton~\cite{eggletonerdosskilton1985colouring} (cf.~\cite[Lem.~2]{barajasserra2009}).


\begin{proposition}\label{prop:primelessk+1}
If $a \le k+1$ is an integer such that $n_j = m_j \, a + r_j$ with $0 < r_j < a$, for $1 \le j \le k$, then $\bn \in \cK(\bm)$.
In particular, if there exists an integer $\le k+1$ that does not divide any of $n_1, n_2, \dots, n_k$, then $\bn$ is a lonely runner instance.
\end{proposition}
\begin{proof}
We need to show that $\cP(\bn) \cap \ZZ^k \ne \emptyset$, which by \eqref{eqnDescriptionPn} means that there exist $\lambda \in [0,1]$ and
$\mu_1, \dots, \mu_k \in [\tfrac{1}{k+1},\tfrac{k}{k+1}]$ such that for all $1 \le j \le k$
\[
  \lambda n_j - \mu_j \in \ZZ \, .
\]
Let $a \le k+1$ be an integer not dividing $n_1, n_2, \dots, n_k$. Then the fractional part $\{ \frac{ n_j }{ a } \}$ lies in $[ \frac 1 a,
\frac{ a-1 }{ a }] \subseteq [\tfrac{1}{k+1},\tfrac{k}{k+1}]$, and so
\[
  \lambda = \frac 1 a \qquad \text{ and } \qquad \mu_j = \left\{ \frac{ n_j }{ a } \right\} 
\] 
will do the trick.
With these choices, we have $m_j = \lambda n_j - \mu_j$ and thus $\bm \in \cP(\bn)$, that is, $\bn \in \cK(\bm)$.
\end{proof}

As a simple consequence of Proposition~\ref{prop:primelessk+1}, for any $\bm \in \ZZ_{ \ge 0 }^k$, we have $2 \, \bm + \bone \in \cK(\bm)$, where $\bone$ is the all-ones vector.
In particular, if all~$n_j$ are odd, then~$\bn$ is a lonely runner instance.
Not much more seems to be known regarding the parities of the speeds.
We generalize the observation just made and, under the given assumptions, we provide an explicit lattice point $\bm$ whose associated cone~
$\cK(\bm)$ contains~$\bn$.

\begin{theorem}\label{thm:oddlrc}
Suppose $n_1 \geq n_2 \geq \dots \geq n_k$, let $E := \{ j \in [k] : \, n_j \text{ \rm is even} \}$, $O := [k] \setminus E$, and
\[
  m_j := \begin{cases}
    \,\,\,\, \frac { n_j } 2 & \text{ if } j \in E, \\ 
    \frac { n_j - 1 } 2 & \text{ if } j \in O.
  \end{cases}
\]
If
\[
  \max \left\{ n_j : \, j \in O \right\} \ \le \ \tfrac{ k-1 }{ 2 } \min \left\{ n_j : \, j \in E \right\}
  \qquad \text{ and } \qquad
  \max \left\{ n_j : \, j \in E \right\} \ \le \ k \, \min \left\{ n_j : \, j \in E \right\} ,
\]
then $\bn \in \cK(\bm)$.
In particular, $\bn$ is a lonely runner instance.
\end{theorem}
\begin{proof}
We claim that $\bn$ satisfies~\eqref{eq:lrcfractions1}.
First, let $i, j \in O$.
Then, by routine manipulations,
\[
  \frac{ (k+1) \, \frac{ n_j - 1 }{ 2 } + 1 }{ (k+1) \, \frac{ n_i - 1 }{ 2 } + k }
  \ = \ \frac{ (k+1) \, n_j - ( k - 1 ) }{ (k+1) \, n_i + ( k - 1 ) } 
  \ \le \ \frac{ n_j }{ n_i }
  \ \le \ \frac{ (k+1) \, n_j + ( k - 1 ) }{ (k+1) \, n_i - ( k - 1 ) } 
  \ = \ \frac{ (k+1) \, \frac{ n_j - 1 }{ 2 } + k }{ (k+1) \, \frac{ n_i - 1 }{ 2 } + 1 }
\]
hold unconditionally.
Second, if $i \in O$ and $j \in E$ then the right inequality in
\[
  \frac{ (k+1) \, \frac{ n_j }{ 2 } + 1 }{ (k+1) \, \frac{ n_i - 1 }{ 2 } + k }
  \ = \ \frac{ (k+1) \, n_j + 2 }{ (k+1) \, n_i + k - 1 } 
  \ \le \ \frac{ n_j }{ n_i }
  \ \le \ \frac{ (k+1) \, n_j + 2k }{ (k+1) \, n_i - k + 1 } 
  \ = \ \frac{ (k+1) \, \frac{ n_j }{ 2 } + k }{ (k+1) \, \frac{ n_i - 1 }{ 2 } + 1 }
\]
holds without conditions, whereas the left inequality requires $n_i \le \frac{ k-1 }{ 2 } \, n_j$.
Finally, if $i, j \in E$ then 
\[
  \frac{ (k+1) \, \frac{ n_j }{ 2 } + 1 }{ (k+1) \, \frac{ n_i }{ 2 } + k }
  \ = \ \frac{ (k+1) \, n_j + 2 }{ (k+1) \, n_i + 2k } 
  \ \le \ \frac{ n_j }{ n_i }
  \ \le \ \frac{ (k+1) \, n_j + 2k }{ (k+1) \, n_i + 2 } 
  \ = \ \frac{ (k+1) \, \frac{ n_j }{ 2 } + k }{ (k+1) \, \frac{ n_i }{ 2 } + 1 }
\]
requires $n_i \le k \, n_j$ and $n_j \le k \, n_i$.
The first condition holds automatically and the second by our assumptions.
\end{proof}

\begin{corollary}
If all but possibly the largest of the $n_j$ are odd, then $\bn$ is a lonely runner instance.
\end{corollary}

\section{Projection Arguments}\label{sect:ProjArgs}

Throughout this part, we assume that $\bn \in \ZZ^k_{>0}$ is such that $n_1 \geq n_2 \geq \dots \geq n_k$.
As we have seen in Proposition~\ref{prop:n1lenk}, the speed vector $\bn$ is a lonely runner instance when $n_1 \leq k \, n_k$.
In the following we show how projections of $\cP(\bn)$ can be used to relax $n_1 \leq k \, n_k$ to conditions that are independent on the fastest runner.
Before we can formulate our result, we describe the polyhedra that arise as projections and cross sections of~$\cP(\bn)$ by coordinate subspaces.

For $\ell \in \{1,\dots,k\}$, let $L_\ell = \lin\{\be_1,\dots,\be_\ell\}$, where $\be_i$ is the $i$th coordinate unit vector.
The orthogonal projection of $\cP(\bn)$ along $L_\ell$ is given by
\begin{align*}
\cP(\bn) \mid L_\ell^\perp \ = \ &\left( \RR \, \bn - \left[\tfrac{1}{k+1},\tfrac{k}{k+1}\right]^k \right) \mid L_\ell^\perp \ = \ \RR \, (n_{\ell+1},\ldots,n_k) - \left[\tfrac{1}{k+1},\tfrac{k}{k+1}\right]^{\{\ell+1,\ldots,k\}} \\
\ = \ &\left\{(x_{\ell+1},\ldots,x_k) : \, \tfrac{1}{k+1} \,  n_i - \tfrac{k}{k+1} \, n_j \ \le \ n_j \, x_i - n_i \, x_j \ \le \ \tfrac{k}{k+1} \, n_i - \tfrac{1}{k+1} \, n_j \ \text{ for }\ell < i < j \leq k \right\}.
\end{align*}
This projection contains the origin if and only if $n_{\ell+1} \leq k \, n_k$.

On the other hand, assuming that $n_{\ell+1} \leq k \, n_k$ and using that the entries of $\bn$ are ordered nonincreasingly,
\begin{alignat}{2}
\cP(\bn) \cap L_\ell \ = \ &\bigg\{ \bx \in \RR^\ell : \ && \tfrac{1}{k+1} \, n_i - \tfrac{k}{k+1} \, n_j \ \le \ n_j \, x_i - n_i \, x_j \ \le \ \tfrac{k}{k+1} \, n_i - \tfrac{1}{k+1} \, n_j \
\text{ for } 1 \le i < j \le \ell , \nonumber \\
&&& n_i - k \, n_j \ \le \ (k+1)\, n_j \, x_i \ \le \ k \, n_i - n_j \
\text{ for } 1 \le i \le \ell < j \le k \bigg\} \nonumber \\ 
\ = \ &\bigg\{ \bx \in \RR^\ell : \ && \tfrac{1}{k+1} \, n_i - \tfrac{k}{k+1} \, n_j \ \le \ n_j \, x_i - n_i \, x_j \ \le \ \tfrac{k}{k+1} \, n_i - \tfrac{1}{k+1} \, n_j \
\text{ for } 1 \le i < j \le \ell , \label{eqn:polySection} \\
&&& n_i - k \, n_k \ \le \ (k+1)\, n_k \, x_i \ \text{ and } \ (k+1)\, n_{\ell+1} \, x_i \ \le \ k \, n_i - n_{\ell+1} \
\text{ for } 1 \le i \le \ell \bigg\} . \nonumber
\end{alignat}

We are now set up to prove our result.
Its second part says that the Lonely Runner Conjecture holds if there is a block of slow runners (with speeds $n_3,\ldots,n_k$) and a block of fast runners (with speeds $n_1,n_2$) such that the fast runners are at least $k$ times faster than the slow runners.

\begin{theorem}\label{thm:mainProjection}
Let $k \geq 3$, let $\bn \in \ZZ_{ >0 }^k$ and assume $n_1 \geq n_2 \geq \dots \geq n_k$.
\begin{enumerate}[{\rm (a)}]
 \item If $n_2 \leq (k-2) \, n_k$, then $\bn$ is a lonely runner instance.
 \item If $n_3 \leq (k-2) \, n_k$ and $n_2 \geq k \, n_3$, then $\bn$ is a lonely runner instance.
\end{enumerate}
\end{theorem}
\begin{proof}
In both cases we aim to ensure the existence of a lattice point in $\cP(\bn)$.

(a): Projecting along the first coordinate direction, we have seen above that $\cP(\bn) \mid \be_1^\perp$ contains the origin if and only if $n_2 \le k \, n_k$.
Therefore, under the assumption $n_2 \le (k - 2) \, n_k$, it suffices to show that the line segment $\cP(\bn) \cap \RR \, \be_1$ contains an integral point.
This holds since its length $L$ is at least one.
Indeed, in view of~\eqref{eqn:polySection}
\begin{align*}
L \ = \ \frac{k}{k+1} \frac{n_1}{n_2} - \frac{1}{k+1} -\frac{1}{k+1} \frac{n_1}{n_k} + \frac{k}{k+1} 
\ = \ \frac{k \, n_1 \, n_k - n_1 \, n_2}{(k+1) \, n_2 \, n_k} + \frac{k-1}{k+1} \ \geq \ \frac{k \, n_k - n_2}{(k+1) \, n_k} + \frac{k-1}{k+1} \ \geq \ 1 .
\end{align*}

(b): Now we project along the first two coordinates.
Assuming that $n_3 \le (k - 2) \, n_k$, the projection $\cP(\bn) \mid L_2^\perp$ contains the origin and the defining inequalities of $\cP(\bn) \cap L_2$ may be labeled by
\begin{align}
\frac{1}{k+1} \, \frac{n_1}{n_k} - \frac{k}{k+1} \ \le \ &\ x_1 \ \le \ \frac{k}{k+1} \,  \frac{n_1}{n_3} - \frac{1}{k+1}\label{P2_x} \\
\frac{1}{k+1} \, \frac{n_2}{n_k} - \frac{k}{k+1} \ \le \ &\ x_2 \ \le \ \frac{k}{k+1} \, \frac{n_2}{n_3} - \frac{1}{k+1}\label{P2_y} \\
\frac{1}{k+1} \, n_1 - \frac{k}{k+1} \, n_2 \ \le \ n_2 \, x_1 &- n_1 \, x_2 \ \le \ \frac{k}{k+1} \, n_1 - \frac{1}{k+1} \, n_2 .
\label{P2_xy}
\end{align}
This defines a symmetric hexagon with shape and $x_2$-coordinates of some of its vertices illustrated in Figure~\ref{fig:HexRef}.
\begin{figure}[ht]
\includegraphics[height=5.5cm]{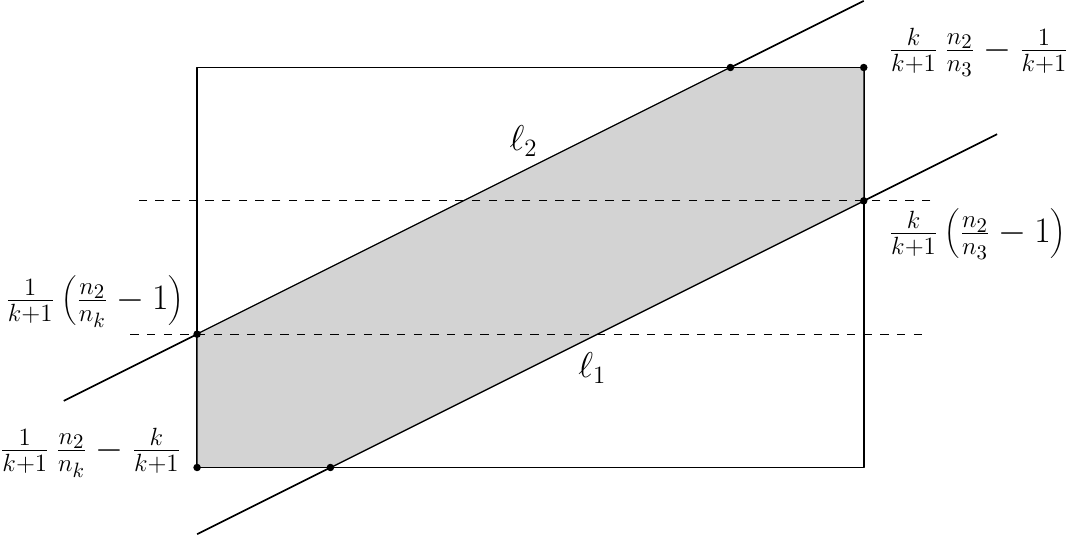}
\caption{The hexagon $\cP(\bn) \cap L_2$ with $x_2$-coordinates of some of its vertices.}\label{fig:HexRef}
\end{figure}
In order to find a lattice point in $\cP(\bn) \cap L_2$ we proceed in two steps.

First, under the assumption $n_2 \ge k \, n_3$, the width of the horizontal strip that is bounded by the dashed lines in Figure~\ref{fig:HexRef} is at least one.
Indeed, assuming also that $n_3 \le (k - 2) \, n_k$ again, we have
\begin{align*}
\frac{k}{k+1} \, \frac{n_2}{n_3} - \frac{k}{k+1} - \frac{1}{k+1} \, \frac{n_2}{n_k} + \frac{1}{k+1} \ = \ \frac{n_2 \left( k \, n_k - n_3 \right) - (k - 1) \, n_3 \, n_k}{(k+1) \, n_3 \, n_k} \ \geq \ \frac{2 \, n_2 \, n_k - (k - 1) \, n_3 \, n_k}{(k+1) \, n_3 \, n_k} \ \ge \ 1 .
\end{align*}
Hence, there exists a horizontal lattice-line that intersects this horizontal strip.

Secondly, we argue that in $x_1$-direction the lines $\ell_1$ and $\ell_2$ are at least of distance one, implying the existence of a lattice point on that very lattice-line.
Said distance $D$ can be computed in view of~\eqref{P2_xy}:
\begin{align*}
D \ = \ \frac{k \, n_1 - n_2 - n_1 + k \, n_2}{(k + 1) \, n_2} \ = \ \frac{k-1}{k+1} \, \frac{n_1 + n_2}{n_2} \ \geq 2 \ \,
\frac{k-1}{k+1} \ \geq \ 1 \, ,
\end{align*}
since $k \geq 3$.

Hence, $\cP(\bn) \cap \ZZ^k \ne \emptyset$ if $n_3 \le (k - 2) \, n_k$ and $n_2 \ge k \, n_3$, so that $\bn$ is indeed a lonely runner instance.
\end{proof}

\section{Cross Section Arguments}\label{sect:crosssectArgs}

A popular line of research is to establish the Lonely Runner Conjecture for speed vectors that are \emph{$L$-lacunary} for a small parameter $L \geq 1$.
Here, a nonincreasing sequence $s_1 \geq s_2 \geq \dots \geq s_k > 0$ is $L$-lacunary, if $\frac{s_j}{s_{j+1}} \geq L$ for $1 \leq i \leq k-1$.
Pandey~\cite{pandey2009anote} started such investigations and proved that if $\bn \in \ZZ^k_{>0}$ is $L$-lacunary with $L = \frac{2 \, (k+1)}{k-1}$, then $\bn$ is a lonely runner instance.
Barajas \& Serra~\cite{barajasserra2009} improved this to $L=2$; Dubickas~\cite{dubickas2011the} achieved $L = 1 + \frac{33 \, \log(k)}{k}$, for large enough $k$; and Czerwi\'{n}ski~\cite{czerwinski2018the} improved slightly on Dubickas' result leaving roughly the $\frac{k+1}{24 \, e}$ slowest runners unconditioned.\footnote{Czerwi\'{n}ski claims $\frac{k+1}{8 \, e}$, but his arguments seem to give only $\frac{k+1}{24 \, e}$.}

An iterative argument based on suitable cross sections of $\cP(\bn)$ yields a proof for a moderate lacunarity of $L = \frac{2 \, k}{k-1}$, but leaving the fastest runner unconditioned.
Notice also that the condition on $\gcd(n_{k-1},n_k)$ below is weaker than the lacunarity condition $\frac{n_{k-1}}{n_k} \ge \frac{2 \, k}{k-1}$.

\begin{theorem}\label{thm:mainCrossSection}
Let $\bn \in \ZZ_{ >0 }^k$ and assume $n_1 \geq n_2 \geq \dots \geq n_k$.
If $\frac{n_{ j }}{n_{ j+1 }} \ge \frac{2 \, k}{k-1} $ for $2 \le j \le k-2$, and $\gcd(n_{ k-1 }, n_k) \le \frac{ k-1 }{ k+1 } (n_{ k-1 } - n_k)$, then $\bn$ is a lonely runner instance.
\end{theorem}
\begin{proof}
We recursively construct a sequence of integers $t_k, t_{ k-1 }, \ldots, t_1$ such that $(t_1, t_2, \ldots, t_k) \in \cP(\bn)\cap\ZZ^k$.
Denote the fractional part of $x \in \RR$ by $\{ x \} := x - \lfloor x \rfloor$.
We can choose $t_k \in \ZZ$ such that
\[
  \left\{ \frac{ n_{ k-1 } }{ n_k } \left( t_k + \frac{ k }{ k+1 } \right) - \frac{ k }{ k+1 } \right\} \ \le \ \frac g { n_k }
\]
where $g := \gcd(n_{ k-1 }, n_k)$.
For example any $t_k \leq -1$ works.
Thus
\[
  t_{ k-1 } \ := \ \left\lfloor \frac{ n_{ k-1 } }{ n_k } \left( t_k + \frac{ k }{ k+1 } \right) - \frac{ k }{ k+1 } \right\rfloor 
\]
satisfies
\[
  \frac{ n_{ k-1 } }{ n_k } \left( t_k + \frac{ k }{ k+1 } \right) - \frac{ k }{ k+1 } -
\frac g { n_k} \ \le \ t_{ k-1 } \ \le \ \frac{ n_{ k-1 } }{ n_k } \left( t_k + \frac{ k }{ k+1 } \right) - \frac{ k }{ k+1 } \, ,
\]
and our condition on $g$ implies that
\begin{equation}\label{eq:doublespeedcondition}
  \frac{ n_{ k-1 } }{ n_k } \left( t_k + \frac{ 1 }{ k+1 } \right) - \frac{ 1 }{ k+1 }
  \ \le t_{ k-1 } \ \le \
  \frac{ n_{ k-1 } }{ n_k } \left( t_k + \frac{ k }{ k+1 } \right) - \frac{ k }{ k+1 } \, .
\end{equation}
Now consider the polytope
\[
  \cP(\bn) \cap \left\{ \bx \in \RR^k : \, x_k = t_k, \ x_{ k-1 } = t_{ k-1 } \right\} 
\]
projected to $\RR^{ k-2 }$, which we call
\[
  \cQ^{ k-2 } \ := \ \left\{ \bx \in \RR^{ k-2 } :
  \begin{array}{l}
    \frac{ n_i - k \, n_j }{ k+1 } \ \le \ n_j \, x_i - n_i \, x_j \ \le \ \frac{ k
\, n_i - n_j }{ k+1 } \ \text{ for } 1 \le i < j \le k-2 \\
    \frac{ n_i }{ n_k } \left( t_k + \frac{ 1 }{ k+1 } \right) - \frac{ k }{ k+1 } \le x_i \le \frac{ n_i }{ n_k } \left( t_k + \frac{ k }{
k+1 } \right) - \frac{ 1 }{ k+1 } \ \text{ for } 1 \le i \le k-2 \\
    \frac{ n_i }{ n_{ k-1 } } \left( t_{ k-1 } + \frac{ 1 }{ k+1 } \right) - \frac{ k }{ k+1 } \le x_i \le \frac{ n_i }{ n_{ k-1 } } \left( t_{ k-1 } + \frac{ k }{ k+1 } \right) - \frac{ 1 }{ k+1 } \ \text{ for } 1 \le i \le k-2
  \end{array}
  \right\} .
\]
By~\eqref{eq:doublespeedcondition},
\[
  \frac{ n_i }{ n_k } \left( t_k + \frac{ 1 }{ k+1 } \right) - \frac{ k }{ k+1 } \ \le \ \frac{ n_i }{ n_{ k-1 } } \left( t_{ k-1 } + \frac{ 1 }{ k+1 } \right) - \frac{ k }{ k+1 }
\]
and
\[
  \frac{ n_i }{ n_{ k-1 } } \left( t_{ k-1 } + \frac{ k }{ k+1 } \right) - \frac{ 1 }{ k+1 } \ \le \ \frac{ n_i }{ n_k } \left( t_k + \frac{ k }{ k+1 } \right) - \frac{ 1 }{ k+1 }
\]
for $1 \le i \le k-2$, and so we can simplify
\[
  \cQ^{ k-2 } \ = \ \left\{ \bx \in \RR^{ k-2 } :
  \begin{array}{l}
    \frac{ n_i - k \, n_j }{ k+1 } \ \le \ n_j \, x_i - n_i \, x_j \ \le \ \frac{ k \, n_i - n_j }{ k+1 } \ \text{ for } 1 \le i < j \le k-2 \\
    \frac{ n_i }{ n_{ k-1 } } \left( t_{ k-1 } + \frac{ 1 }{ k+1 } \right) - \frac{ k }{ k+1 } \le x_i \le \frac{ n_i }{ n_{ k-1 } } \left( t_{ k-1 } + \frac{ k }{ k+1 } \right) - \frac{ 1 }{ k+1 } \ \text{ for } 1 \le i \le k-2
  \end{array}
  \right\}
\]
and revise our goal to prove that $\cQ^{ k-2 } \cap \ZZ^{ k-2 } \ne \emptyset$.

By our assumption that $\frac{n_{ k-2 }}{n_{ k-1 }} \ge \frac{2 \, k}{k-1}$, there exists an integer $t_{ k-2 }$ that satisfies
\[
  \frac{ n_{ k-2 } }{ n_{ k-1 } } \left( t_{ k-1 } + \frac{ 1 }{ k+1 } \right) - \frac{ 1 }{ k+1 } \ \le t_{ k-2 } \ \le \
  \frac{ n_{ k-2 } }{ n_{ k-1 } } \left( t_{ k-1 } + \frac{ k }{ k+1 } \right) - \frac{ k }{ k+1 } ,
\]
since this interval is of length at least one.
So, we can repeat the construction to obtain polytopes
\[
  \cQ^{ k-3 }, \cQ^{ k-4 }, \ldots, \cQ^{ 1 } \ = \ \left\{ x_1 \in \RR : \frac{ n_1 - k \, n_2 }{ k+1 } \ \le \ n_2 \, x_1 - n_1 \, t_2 \ \le \ \frac{ k \, n_1 - n_2 }{ k+1 } \right\} .
\]
The latter is an interval with length $\frac{ k-1 }{ k+1 } (n_1 + n_2) \ge 1$ and thus contains an integer~$t_1$.
\end{proof}

\section{Musings}\label{sec:musings}

The coordinates of an integral point $\bm \in \cP(\bn) \cap \ZZ^k_{\geq 0}$ have the following meaning:
There is a time at which all runners are at least $\frac{1}{k+1}$ away from the starting point and the $i$th runner (with speed
$n_i$) is in her $(m_i+1)$st round on the track.
It is well known that for $k=2$, it happens during the first round of the slower runner that the distance of both runners from the start is at least $\frac13$.
Thus, the cones $\cK(\bm)$ for $\bm \in \ZZ^2_{\geq 0}$ with $m_1 \, m_2 = 0$ already cover the whole nonnegative orthant $\RR^2_{\geq 0}$.
Another setting where this phenomenon occurs is when the slowest runner runs with speed $1$.
Indeed, a result of Czerwi\'{n}ski \& Grytczuk~\cite{czerwinskigrytczuk2008invisible} says that the maximal distance from the starting point that all runners achieve simultaneously is attained at a time $t = \frac{a}{n_i + n_j}$, for some $1 \leq i < j \leq k$ and $a \in \{1,\ldots,n_i + n_j - 1\}$, hence during the first round of the slowest runner.
Aside from these particular cases, we do not know what happens in general.

\begin{question}
Assume that the Lonely Runner Conjecture holds in dimension $k$, that is, $\RR^k_{\geq 0} = \bigcup_{\bm \in \ZZ^k_{\geq 0}} \cK(\bm)$.
Is it true that $\RR^k_{\geq 0}$ is covered by the cones $\cK(\bm)$, where $\bm \in \ZZ^k_{\geq 0}$ runs over the integer points such that $m_i \leq c_k$, for some $1 \leq i \leq k$ and some constant $c_k$ only depending on $k$?
Can $c_k$ be chosen to be $0$?
\end{question}

During our studies of the lonely runner polyhedron $\cP(\bn)$ the following conjecture emerged.
It claims that not only $\cP(\bn)$ but \emph{each of its translates} contains a lattice point, provided that the speeds listed in~$\bn$ are pairwise distinct.
If it could be shown to be equivalent to the Lonely Runner Conjecture, it would mean that the assumption that the runners all start at the same place is unnecessary.
This in fact was recently conjectured by J\"org M.~Wills (personal communication).
The geometric argument for two runners at the beginning of Section~\ref{sect:geometricProofs} shows that the claim holds for $k=2$.

\begin{conjecture}\label{conj:allTranslates}
Let $\bn \in \ZZ_{ >0 }^k$ be such that $n_i \neq n_j$, for every $i \ne j$.
Then, for every translation vector $\bt \in \RR^k$, we have $\left( \cP(\bn) + \bt \right) \cap \ZZ^k \ne \emptyset$.
\end{conjecture}

In view of the equivalences in Proposition~\ref{prop:geometricLRC} the validity of this conjecture would also mean that the translates of the zonotope $\cZ(\bn)$ by vectors of the projected lattice $\ZZ^k \mid \bn^\perp$ cover the hyperplane $\bn^\perp$.
In other words, the covering radius (see, e.g.,~\cite[Ch.~23]{gruber}) of~$\cZ(\bn)$ with respect to $\ZZ^k \mid \bn^\perp$ is
bounded above by one.
Note that the assumption that the speeds $n_i$ are pairwise distinct is crucial.
In fact, the statement of Conjecture~\ref{conj:allTranslates} is not valid for $\bn = (1,1,1)$, for instance.

If we relax both assumptions, that is, we allow the runners to start at different places \emph{and} have equal speeds, then the problem does change.
Even more, the resulting question has been answered already in 1976 in work by Schoenberg~\cite{schoenberg1976extremum}.
He proved that in this setting we need to change the \emph{gap of loneliness} from~$\frac{1}{k+1}$ to the smaller value $\frac{1}{2 \, k}$, and that this is tight.
Of course this implies the Lonely Runner Conjecture for the same bound~$\frac{1}{2 \, k}$, which had been shown much earlier by Wills~\cite{willslonelyrunner}.
It seems that it has not been noticed that Wills' application of the union bound, in turn, also gives a slick proof of Schoenberg's result.

\begin{theorem}[Schoenberg 1976]
Given positive integers $n_1, n_2, \dots, n_k$ and reals $s_1, s_2, \dots, s_k$, there exists a real number $t$ such that for all $1 \le j \le k$ the distance of $s_j + t \, n_j$ to the nearest integer is at least $\frac{ 1 }{ 2 \, k }$.
Furthermore, this bound cannot be improved for $n_i = 1$ and $s_i = \frac{i-1}{k}$, for $1 \leq i \leq k$.
\end{theorem}
\begin{proof}
Let $\lambda \in [0,\frac12]$.
The distance of $s_j + t \, n_j$ to the nearest integer is at least $\lambda$ if and only if $s_j + t \, n_j \in \ZZ + [\lambda, 1 - \lambda]$.
By the periodicity of the problem it suffices to look at $t \in [0,1]$.
Define
\[
  I_j \ := \ [s_j, s_j + n_j] \cap (\ZZ + [\lambda, 1 - \lambda]) - s_j \, ,
\]
which is a union of closed intervals.
The crucial observation is that the total length of~$I_j$ is independent of $s_j$.
Indeed,
\[
\left| I_j \right| \ = \ \left| [s_j, s_j + n_j] \cap (\ZZ + [\lambda, 1 - \lambda]) \right| \ = \ \left| [0, n_j] \cap (\ZZ +
[\lambda, 1 - \lambda]) \right| \ = \ n_j \, (1 - 2 \, \lambda) \, ,
\]
because the $n_j$ are integral.
Now, the union bound in elementary probability theory implies that
\begin{align*}
\PP\left( \bigcup_{j=1}^k \left\{ t \in [0,1] : t \, n_j \notin I_j \right\} \right) \ \leq \ \sum_{j=1}^k \PP\left( t \in [0,1] : t \, n_j \notin I_j \right) \ = \ \sum_{j=1}^k \left( 1 - \frac{\left| I_j \right|}{n_j} \right) \ = \ 2 \, k \, \lambda .
\end{align*}
Hence, there is a desired real number $t \in [0,1]$ whenever $\lambda < \frac{1}{2 \, k}$.
By the compactness of the $I_j$ this is also true for $\lambda = \frac{1}{2 \, k}$.
\end{proof}


\bibliographystyle{amsplain}
\bibliography{bib}

\setlength{\parskip}{0cm} 

\end{document}